\documentclass[11pt]{amsart}

\newtheorem{thm}{Theorem}[section]
\newtheorem{lem}[thm]{Lemma}

\newtheorem{cor}[thm]{Corollary}
\newtheorem{prop}[thm]{Proposition}

\theoremstyle{plain} 
\newtheorem{main}{Theorem}
\newtheorem{corm}[main]{Corollary}

\theoremstyle{definition} 

\newtheorem{defn}[thm]{Definition}
\newtheorem{rem}[thm]{Remark}

\newtheorem*{remar}{Remark}

\numberwithin{equation}{section}

\usepackage{amssymb,latexsym,amscd,amsmath,amsthm}
\usepackage[all]{xy}
\usepackage{epsfig}

\newcommand{\R}{\mathbf{R}}
\newcommand{\Z}{{\mathbf{Z}}}
\newcommand{\Q}{\mathbf{Q}}
\newcommand{\C}{{\mathbf{C}}}

\newcommand{\sph}{\mathbf{S}}
\newcommand{\f}{\mathbf{F}}

\newcommand{\SU}{\operatorname{SU}}

\newcommand{\syp}{\operatorname{Sp}}
\newcommand{\SO}{\operatorname{SO}}

\newcommand{\U}{\operatorname{U}}
\newcommand{\Hom}{\operatorname{Hom}}

\newcommand{\id}{\operatorname{{\it id}}}

\newcommand{\tr}{\operatorname{tr}}

\renewcommand{\lim}[1]{\mathop{\underset{#1} {\underset \longleftarrow
{\text{\rm lim}}}}}

\newcommand{\ad}{\operatorname{ad}}

\newcommand{\I}{{\bf I}}

\newcommand{\ek}{{\rm ek}}
\newcommand{\wu}{{\rm Wu}}

\newcommand{\LC}{^{\rm LC}}
\newcommand{\G}{{{\rm G}_2}}
\newcommand{\Spin}{{\rm Spin}}
\newcommand{\Cl}{{\rm Cl}}
\newcommand{\End}{{\rm End}}
\newcommand{\even}{{\rm ev}}
\newcommand{\Ca}{\mathbf{O}}
\newcommand{\so}{\mathfrak{so}}
\newcommand{\into}{\hookrightarrow}
\newcommand{\frh}{\mathfrak{h}}
\newcommand{\frg}{\mathfrak{g}}
\newcommand{\frp}{\mathfrak{p}}
\newcommand{\frs}{\mathfrak{s}}
\newcommand{\frt}{\mathfrak{t}}
\newcommand{\fru}{\mathfrak{u}}
\newcommand{\g}{\mathfrak{g}_2}
\newcommand{\Vol}{\operatorname{vol}}
\newcommand{\<}{\langle}
\renewcommand{\>}{\rangle}
\newcommand{\phy}{{\varphi}}
\newcommand{\Adach}{{\widehat{A}}}
\newcommand{\Ldach}{{\widehat{L}}}
\newcommand{\adach}{{\hat{A}}}
\newcommand{\ldach}{{\hat{L}}}
\newcommand{\Adachsl}{{\,\,\widetilde{\!\!\widehat{A}}}}
\newcommand{\Ldachsl}{{\,\widetilde{\!\widehat{L}}}}
\newcommand{\adachsl}{{\tilde{\hat A}}}
\newcommand{\ldachsl}{{\tilde{\hat L}}}
\newcommand{\ch}{{\operatorname{ch}}}
\newcommand{\punkt}{\,\cdot\,}
\newcommand{\norm}[1]{\left\|#1\right\|}
\newcommand{\abs}[1]{|#1|}
\newcommand{\Cal}[1]{\mathcal{#1}}
\newcommand{\pix}{{\tilde{\pi}}}
\newcommand{\piy}{{\hat{\tilde{\pi}}}}
\newcommand{\aD}{{\,\widetilde{\!{\rm ad}\!}\,}}
\newcommand{\aDsl}{{\,\widehat{\widetilde{\!{\rm ad}\!}}\,}}
\newcommand{\adX}{{\aD_\frp}}
\newcommand{\adY}{{\aDsl_\frp}}
\newcommand{\Dsl}{\widetilde{D}}
\newcommand{\Bsl}{\widetilde{B}}

\begin{document}

\newcommand{\spacing}[1]{\renewcommand{\baselinestretch}{#1}\large\normalsize}
\spacing{1.14}

\title[Diffeomorphism type of the Berger space]{Diffeomorphism type of
the Berger space $\SO(5)/\SO(3)$}

\author[S.\ Goette, N.\ Kitchloo \& K.\ Shankar]{Sebastian
Goette$^\ast$, Nitu Kitchloo, Krishnan Shankar$^{\ast\ast}$}

\thanks{$^\ast$ Research at MSRI supported in part by NSF grant
DMS-9701755.}
\thanks{$^{\ast \ast}$ Research supported in part by NSF grant
DMS-0103993.}

\address{Mathematisches Institut der Uni T\"ubingen,
Auf der Morgenstelle~10, 72076 T\"ubingen, Germany.}
\email{sebastian.goette@uni-tuebingen.de}
\address{Department of Mathematics\\ Johns Hopkins University\\
Baltimore\\ MD 21218.}
\email{nitu@math.jhu.edu}
\address{Department of Mathematics\\ University of Michigan\\ Ann
Arbor\\ MI 48109.}
\email{shankar@umich.edu}

\begin{abstract} We compute the Eells-Kuiper invariant of the
Berger manifold~$\SO(5)/\SO(3)$ and determine that it is diffeomorphic
to the total space of an $\sph^3$-bundle over~$\sph^4$. This answers a
question raised by K.\ Grove and W.\ Ziller.
\end{abstract}

\maketitle

\normalsize
\thispagestyle{empty}

\section*{Introduction}

There has been renewed interest in Riemannian structures of
non-negative or positive curvature on the total spaces of
$\sph^3$-bundles over~$\sph^4$. These bundles have been of interest to
topologists since Milnor showed that if the Euler class of such a
bundle is~$\pm 1$, then the total space is a homotopy sphere. Until
recently there was only one exotic sphere, the so called
Gromoll-Meyer sphere (cf.\ \cite{gm}), which was known to admit a
metric of non-negative curvature. Then in their paper \cite{gz}, K.\
Grove and W.\ Ziller showed that every $\sph^3$-bundle over~$\sph^4$
admits infinitely many complete metrics of non-negative curvature. In
particular, all the exotic Milnor spheres admit such metrics. Which of
course begs the question: which exotic spheres, or more generally,
which $\sph^3$-bundles over~$\sph^4$ admit metrics of positive
sectional curvature?

The Berger space, $M^7 = \SO(5)/\SO(3)$, was first described by M.\
Berger as a manifold that admits a (normal) homogeneous metric of
positive sectional curvature. The embedding of~$\SO(3)$ in~$\SO(5)$ is
maximal and irreducible (cf.\ \cite{wolf}), it is a rational homology
sphere with~$H^4(M,\Z) = \Z_{10}$ (cf.\ \cite{berger}) and it has the
cohomology ring of an $\sph^3$-bundle over~$\sph^4$. In \cite{gz}, K.\
Grove and W.\ Ziller asked whether the Berger space is topologically
or differentially equivalent (as a manifold) to an $\sph^3$-bundle
over~$\sph^4$. Part of this was settled in \cite{kish} where it was
shown that the Berger space is $PL$-homeomorphic to such a bundle. To
settle the diffeomorphism question requires computing the
Eells--Kuiper invariant~$\ek(M)$.  The original definition of~$\ek(M)$
in~\cite{ek} requires that $M^7$ be written as the boundary of an
eight dimensional spin manifold. Since the cobordism group~$\Omega_7$
is known to be trivial, one knows that any closed, 2-connected,
7-manifold admits a spin coboundary. However, an explicit coboundary
for the Berger space has not been found.

Instead, we use the analytic formula~\eqref{EK} for the Eells--Kuiper
invariant due to Donnelly~\cite{donnelly} and Kreck--Stolz~\cite{ks},
which is based on the Atiyah--Patodi--Singer index theorem~\cite{aps}.
This formula expresses~$\ek(M)$ in terms of the $\eta$-invariants of
the signature operator and the untwisted Dirac operator.  To determine
these $\eta$-invariants, we follow the first named author's approach
in~\cite{g1}, \cite{g2}, \cite{g3}; we first replace the operators~$B$
and~$D$ that are associated to the Levi-Civita connection on~$M$ by
operators of the kind~$D^{\frac13}$.  These operators are particularly
well adapted to representation theoretic methods, which was first
noticed by Slebarski in \cite{sleb}, and later exploited in \cite{g1},
\cite{g2} and \cite{kostant}. In particular, we use the explicit
formula for $\eta$-invariants from \cite{g3}.  As a consequence of
this approach, we can employ the reductive connection on~$M$ to
compute the secondary Pontrjagin number in \eqref{EK} which is easier
than using the Levi-Civita connection. This is accomplished by
exploiting the existence of an equivariant~${\rm G}_2$ structure on
the tangent bundle~$TM$. Our main result is:

\begin{main}\label{MainThm}
Depending on the orientation the Eells-Kuiper invariant
(EK-invariant) of the Berger space~$M=\SO(5)/\SO(3)$ is
	$$\ek(M)=\pm\frac{27}{1120}.$$
\end{main}

Note that for our choice of orientation (which descends from a choice
of orientation on the octonions) the exact value we obtain is $\ek(M)
= -\frac{27}{1120}$. By \cite{kish} we know that the Berger space may
be given an orientation so as to have the oriented $PL$-type of some
$\sph^3$-bundle~$M_{m,10}$ over~$\sph^4$ with Euler class~$10$ and
Pontrijagin class~$2\,(10+2m)\in\Z\cong H^4(\sph^4)$ with respect to
the standard generator (in the notation of~\cite{ce}).  In fact,
because the value of the $PL$-invariant~$s_1(M)=28\,\ek(M)\in\Q/\Z$
of~\cite{ks} equals~$\frac{13}{40}$, it follows from Theorem~1.2
of~\cite{ce} by an explicit computation of all possible values
of~$s_1(M_{m,10})$ that the Berger space~$M$ is orientation preserving
(reversing) $PL$-equivalent to~$M_{m,10}$ if and only if $m\equiv\pm2$
($m\equiv\pm1$) modulo~$10$.

Given a pair of 2-connected, 7-manifolds $M_1$ and $M_2$, they are
$PL$-homeomorphic to each other if and only if there exists an exotic
sphere $\Sigma$ so that $M_1 \# \Sigma = M_2$. This is a consequence
of the fact that $PL/O$ is 6-connected (cf.\ \cite{mm}).  Moreover,
the $EK$-invariant is additive with respect to connected sums and
attains 28 distinct values on the group of exotic 7-spheres. The
previous two facts were used in \cite{ce} to do the diffeomorphism
classification of $\sph^3$-bundles over $\sph^4$.
Hence, $M_1$ and $M_2$ are oriented diffeomorphic if and only if they
are $PL$-homeomorphic and have the same $EK$-invariant.
Comparing this with the values of~$\ek(M_{m,n})$
in~\cite{ce}, we obtain the following corollary.

\begin{corm}\label{cor2} The Berger space is diffeomorphic
to~$M_{\mp1,\pm10}$, the $\sph^3$-bundle over~$\sph^4$ with Euler
class~$\pm 10$ and first Pontrjagin class equal to $\pm16$ times the
generator in~$H^4(\sph^4)$ with respect to the standard choice of
orientation on~$\sph^4$. \end{corm}

\begin{remar}
In general, any $\sph^3$-bundle over~$\sph^4$ with non-vanishing Euler
class~$n\in\Z\cong H^4(\sph^4)$ is diffeomorphic to infinitely many
other $\sph^3$-bundles over $\sph^4$ with the same Euler class.  It
follows from Corollary~1.6 in~\cite{ce} that the Berger space with the
orientation specified in~\eqref{2.4} is orientation reversing
diffeomorphic to $M_{m,n}$ (orientation preserving diffeomorphic
to~$M_{-m,-n}$) if and only if $n=10$ and $m$
is congruent modulo $140$ to $-1$, $-9$, $-29$ or $19$;
this was pointed out to us by C.\ Escher.
Note that there is no space $M_{m,10}$ that is orientation reversing
diffeomorphic to $M_{-1,10}$.
\end{remar}

We also mention another consequence of Theorem~\ref{MainThm}. It is a
natural question to ask: what is the largest \textit{degree of
symmetry} for $\sph^3$-bundles over~$\sph^4$? The degree of symmetry
of a Riemannian manifold is the dimension of its isometry group. For
instance, it is well known that the maximal degree of symmetry for
exotic 7-spheres is 4 (cf.\ \cite{straume}). However, some
$\sph^3$-bundles over~$\sph^4$ admit actions of larger groups.  It
follows from~\cite[Theorem 4]{onishchik} (see also \cite{klaus}), that
apart from the trivial bundle~$\sph^4\times\sph^3$, the only seven
dimensional homogeneous manifolds that have the cohomology of an
$\sph^3$-bundle over~$\sph^4$ are~$\sph^7 = \SO(8)/\SO(7)$, $T_1\sph^4
= \syp(2)/\Delta\syp(1)$, the unit tangent bundle of~$\sph^4$, the
Berger space~$M = \SO(5)/\SO(3)$, as observed by K.\ Grove and W.\
Ziller in \cite{gz}.
The spaces~$\sph^4 \times \sph^3$, $\sph^7$ and $T_1 \sph^4$ are
diffeomorphic to principal $\sph^3$-bundles over~$\sph^4$.  On the
other hand, it was shown in \cite{gz} that the Berger space is not
diffeomorphic (or even homeomorphic) to a principal $\sph^3$-bundle
over $\sph^4$ (since its first Pontrjagin class does not vanish), but
it is homotopy equivalent to a principal $\sph^3$-bundle
over~$\sph^4$.

\begin{corm}\label{Cor3}
Up to diffeomorphism,
the only total spaces of $\sph^3$-bundles over~$\sph^4$
that are homogeneous are the trivial
bundle, the Hopf bundle, the unit tangent bundle
of~$\sph^4$, and the Berger space~$\SO(5)/\SO(3)$. \end{corm}

The paper is organized as follows: in Section \ref{s1} we motivate the
question of the diffeomorphism type from problems in the geometry of
positive curvature. In Section \ref{s2} we compute the EK-invariant
for the Berger space using spectral theory. In Section \ref{s3} we
discuss the existence of independent vector fields on 2-connected
7-manifolds.

The first named author would like to thank Marc v.\ Leeuwen for his generous
help with the computer program LiE \cite{lie}. We would also like to
thank Wolfgang Ziller for useful comments.

\section{Motivation}\label{s1}

It is a general problem in Riemannian geometry to find and describe
closed manifolds that admit a metric of positive sectional
curvature. There are few known obstructions and frustratingly, few
known examples. The difficulty stems from the fact that all known
examples arise as quotients of Lie groups --- as homogeneous spaces or
as biquotients (double coset manifolds). Simply connected, homogeneous
manifolds with positive curvature were classified by Berger, Wallach
and Berard-Bergery in the sixties and seventies. The Berger space
evidently appears in the classification of normal homogeneous
manifolds of positive curvature due to M.\ Berger \cite{berger}. Other
than the homogeneous spaces of positive curvature there are some
examples in low dimensions, but in dimensions 25 and up the only known
examples are the compact, rank one, symmetric spaces.

The only way we know to construct examples of positively curved
manifolds is to look at quotients of compact Lie groups. By the
Gray-O'Neill curvature formulas, submersions are curvature
non-decreasing. So one looks for positive curvature at the base of a
Riemannian submersion. However, all known examples of positively
curved manifolds, except the Berger space, fit into a fibration
sequence, like the Hopf fibration of spheres over projective
spaces. Fibrations may provide us with another way to construct
examples of positively curved manifolds by the following method: Given
a principal $G$-bundle, $G\rightarrow P\rightarrow B$, a
\textit{connection metric} on~$P$ is a choice of principal connection
$\omega$ i.e., a choice of horizontal space~$\mathcal{H}_G$ in~$P$
invariant under~$G$ such that the map~$P\rightarrow B$ is a Riemannian
submersion with totally geodesic fibers. The fibers are all isometric
to each other and the metric on any fiber is isometric to~$(G,\langle
, \rangle)$ for some choice of left invariant metric on~$G$.
By Hermann~\cite{hermann},
every submersion metric on~$P$ with totally geodesic fibers
must be of this form. Now we look at
associated bundles~$G/H \rightarrow M=P\times_G G/H \rightarrow B$
with fiber~$G/H$. We declare the fibers to be orthogonal to the
horizontal spaces~$\mathcal{H}$, where~$\mathcal{H}$ in~$M$ is the
image of~$\mathcal{H}_g \times \{0\} \subset T(P\times G/H)$. The
metric on the total space is taken to be the orthogonal sum of the
metrics on the fibers, where each fiber is isometric to~$G/H$ with a
normal homogeneous metric (or more generally a left invariant metric),
and the pullback of the metric on the base.

If we have a fibration with a connection metric, then the fiber~$G/H$,
which is totally geodesic, must be a circle or a normal homogeneous
space of positive curvature. All known homogeneous spaces of positive
curvature fit into fibrations with connection metrics, except the
Berger space. Furthermore, Derdzinski and Rigas have shown in
\cite{dr} that for $\sph^3$-bundles over~$\sph^4$, the only bundle
that admits a connection metric of positive curvature is the Hopf
bundle whose total space is the round sphere. Since we now know that
the Berger space is diffeomorphic to the total space of an
$\sph^3$-bundle over~$\sph^4$, it follows that its metric is not a
connection metric. If one could find an explicit smooth submersion to
$\sph^4$, then we could check whether the positive curvature metric is
a submersion metric. At the very least it makes plausible the
suggestion that there are more general metrics of positive curvature
on bundles with large degree of symmetry that are not connection
metrics.

\section{The Eells-Kuiper invariant of $\SO(5)/\SO(3)$}\label{s2}

To compute the Eells-Kuiper invariant we use the formula
\begin{equation}\label{EK}
\ek(M)
	=\frac{\eta(B)}{2^5\,7}+\frac{\eta(D)+h(D)}2
		-\frac{1}{2^7\,7}\,\int_Mp_1(M,\nabla\LC)\wedge h(M,\nabla\LC)
	\qquad\in\Q/\Z
\end{equation}
due to Donnelly \cite{donnelly} and Kreck-Stolz \cite{ks}. Here~$B$
and~$D$ are the odd signature operator and the untwisted Dirac
operator on~$M$, and~$h(M,\nabla\LC)\in\Omega^3(M)$ is a form whose
exterior differential is the first Pontrjagin form~$p_1(M,\nabla\LC)$
with respect to the Levi-Civita connection.  Equation \eqref{EK} has the
advantage that we do not need to find an explicit zero spin bordism
for~$M$.

Now we perform all the computations necessary to determine a numerical
value for the Eells-Kuiper invariant~$\ek(M)$ for~$M=\SO(5)/\SO(3)$
using the methods of \cite{g1}, \cite{g2}, \cite{g3}.  In Section
\ref{s2.1} we recall the $\G$-structure on~$TM$, which will be
important for calculations throughout this chapter. In Section
\ref{s2.2} we control the spectral flow of the deformation of the odd
signature operator to Slebarski's ${\frac13}$-operator. In Section
\ref{s2.3} we determine the $\eta$-invariants of the Dirac operator
and the odd signature operator on~$M$ up to a local correction. In
Section \ref{s2.4} we adapt~\eqref{EK} to our situation.  Finally in
Section \ref{s2.5} we compute the remaining local correction term and
obtain the value of~$\ek(M)$.

\subsection{The $\G$-structure on $TM$}\label{s2.1}\hfill

Using Schur's lemma, we exhibit a $\G$-structure on the tangent bundle
of~$M=\SO(5)/\SO(3)$.  Using this structure, we will be able to
simplify several explicit calculations needed to control both the
equivariant spectral flow from the Riemannian signature operator~$B$
to its reductive (or ``cubical'') deformation~$\widetilde B$, and the
Chern-Simons correction term.  We will also use some branching rules
for~$\SO(3)\subset\G$; these can be checked using a suitable computer
program like LiE \cite{lie}.

To facilitate computations, let~$e_{ij}\in\so(n)$ for~$i\ne j$ denote
the endomorphism that maps the $j$-th vector~$e_j$ of the standard
orthonormal base of~$\R^n$ to the $i$-th vector~$e_i$ and~$e_i$ to
$-e_j$, and vanishes on all other vectors.  Then~$e_{ji}=-e_{ij}$ and
$[e_{ij},e_{jk}]=e_{jk}$ unless~$i=k$.  We fix the scalar product
$\<A,B\>=-{\frac12}\,\tr(AB)$ on~$\mathfrak{so}(5)$, so that the basis
$e_{ij}$ of~$\so(5)$ becomes orthonormal.  We fix an embedding
$\iota\colon\so(3)\into\so(5)$ with
\begin{equation}\label{2.2}
\begin{aligned}
	e_{12}\mapsto\iota_{12}
	&=2e_{12}+e_{34}\;,\\
	e_{23}\mapsto\iota_{23}
	&=e_{23}-e_{14}+\sqrt3\,e_{45}\;,\\ \text{and }\quad
	e_{13}\mapsto\iota_{13}
	&=e_{13}+e_{24}+\sqrt3\,e_{35}\;.
\end{aligned}
\end{equation}
Then the vectors $\iota_{ij}$ of $\frh\cong\so(3)$ are orthogonal and
of length $\sqrt 5$.  Let $\frp$ be the orthogonal complement of
$\frh=\iota(\so(3))$ in $\frg=\so(5)$.

The embedding of $H=\SO(3)$ in $\SO(5)$ for the Berger space is given
by the conjugation action of $\SO(3)$ on real, $3\times 3$ symmetric
matrices of trace zero. The isotropy representation $\pi$ of $H$ on
$\frp$ is the seven-dimensional irreducible representation of $\SO(3)$
and so the Berger space is isotropy irreducible (cf.\ \cite{wolf}).  It
is well known that the seven-dimensional, irreducible, orthogonal
representation of $\SO(3)$ factors through the groups $\G$ and
$\Spin(7)$,
$$
	\SO(3)\to\G\to\Spin(7)\to\SO(7)\;.
$$
The second factorization is due to the fact that $\G$ is simply
connected; it implies that $M$ admits an $\SO(5)$-equivariant spin
structure. The first factorization is more important, since it allows
us to identify $\frp$ with the space of imaginary octonions $\I$; here
we write Cayley's octonions $\Ca=\R\oplus\I$ as split into their real
and imaginary parts.  Let $*$ denote the Cayley product, and let
$*_\I$ denote its projection onto $\I$.  Note that $\Ca$ carries a
natural scalar product $\<p,q\>$ given by the real part of $p\overline
q$, and that $\G$ and $H$ preserve the decomposition $\Ca=\R\oplus\I$
as well as $*$ and~$\<\punkt,\punkt\>$.

\begin{lem}\label{l2.1}
With a suitable isometric, $\G$-equivariant identification
of~$\frp$ with the imaginary octonions,
one has
	$$[v,w]_\frp={\frac1{\sqrt5}}\,v*_\I w
		\text{ for all~$v$, $w\in\frp$.}$$
\end{lem}
\begin{proof} Schur's lemma implies that $[v,w]_\frp=c\,v*_\I w$
for some real constant $c$, because $\Lambda^2\frp$ splits
$\SO(3)$-equivariantly into the irreducible real
$\SO(3)$-representations $\kappa_1$, $\kappa_3$ and $\kappa_5$ of
dimensions $3$, $7$ and~$11$, each of multiplicity one.  On the other
hand, $\kappa_3$ is just the restriction of the standard
representation of $\G$, which leaves ``$*_\I$'' invariant.

To determine $c$, we pick two orthogonal unit vectors $v$, $w\in\frp$.
Because then $v*w\in\I$, we have $\norm{v*w}=1$ and
$\norm{[v,w]_\frp}=\abs c$.  For example with $v={\frac1{\sqrt5}}\,
(e_{12}-2e_{34})$ and $w=e_{25}\in\frp$, we find
$[v,w]={\frac1{\sqrt5}}\,e_{15}\in\frp$, so
$$
	\pm c=\norm{[v,w]_\frp}={\frac1{\sqrt5}}\;.
$$
Because $\frp$ is irreducible, an isometric, $\G$-equivariant
identification $\I\cong\frp$ is unique up to sign, and we may pick the
sign so that $c={\frac1{\sqrt5}}$.
\end{proof}

For later use, we explicitly identify $\frp\cong\I$ as in Lemma~\ref{l2.1}.
We are also interested in the decomposition of $\frp\otimes\C$ into
weight spaces for the $H$-representation $\pi$.  Recall that $\I$
admits an orthonormal base $e_1$, \dots, $e_7$ such that
\begin{equation}\label{2.3}
	e_i*e_{i+1}=e_{i+3}
\end{equation}
for all $i\in\{1,\dots,7\}$, where the indices $i+1$ and $i+3$ are to
be understood modulo $7$.  We identify these imaginary octonions with
$\frp$ with an orthonormal basis given by
\begin{equation}\label{2.4}
\begin{aligned}
	e_1={\frac1{\sqrt5}}\, e_{12}-\frac2{\sqrt5}\,e_{34}\;,\qquad
	&e_2=\frac{\sqrt2}{\sqrt5}\,e_{45}
		-\frac{\sqrt3}{\sqrt{10}}\,(e_{23}-e_{14})\;,\qquad
	e_3=e_{25}\;,\\
	e_4=\frac{\sqrt2}{\sqrt5}\,e_{35}
		-\frac{\sqrt3}{\sqrt{10}}\,(&e_{13}+e_{24})\;,\qquad
	e_5=\frac1{\sqrt2}\,(e_{24}-e_{13})\;,\\
	e_6=-\frac1{\sqrt2}\,(&e_{23}+e_{14})\;,
		\quad\text{and}\quad
	e_7=e_{15}.
\end{aligned}
\end{equation}
We leave it to the reader to check that indeed $[e_i,e_{i+1}] =
{\frac1{\sqrt5}}\,e_{i+3}$ where the indices are taken modulo $7$.

Let us also compute the action of $\frh$ on $\frp$.  If $\iota_{12}$,
$\iota_{13}$ and $\iota_{23}$ are given by \eqref{EK},
then $f_1={\frac1{\sqrt5}}\,\iota_{12}$,
$f_2={\frac1{\sqrt5}}\,\iota_{23}$ and $f_3={\frac1{\sqrt5}}\,\iota_{13}$
form an orthonormal base of $\frh$.  For~$k=1$, \dots, $3$, we define
an element of $\Lambda^2TM$ by
$$
\alpha_k=\<f_k,[\punkt,\punkt]\>=\<\pi_{*f_k}(\punkt),\punkt\>,
$$
so
\begin{equation}\label{2.5}
\begin{aligned}
	\alpha_1
	&={\frac1{\sqrt5}}\,e^2\wedge e^4+\frac2{\sqrt5}\,e^3\wedge e^7
		-\frac3{\sqrt5}\,e^5\wedge e^6\;,\\
	\alpha_2
	&=\frac{\sqrt6}{\sqrt5}\,e^1\wedge e^4-\frac1{\sqrt2}\,e^2\wedge e^7
		-\frac1{\sqrt2}\,e^3\wedge e^4
		+\frac{\sqrt3}{\sqrt{10}}\,e^3\wedge e^5
		+\frac{\sqrt3}{\sqrt{10}}\,e^6\wedge e^7\;,\\
			\text{and }\quad
	\alpha_3
	&=\frac{\sqrt6}{\sqrt5}\,e^1\wedge e^2-\frac1{\sqrt2}\,e^2\wedge e^3
		+\frac{\sqrt3}{\sqrt{10}}\,e^3\wedge e^6
		-\frac1{\sqrt2}\,e^4\wedge e^7
		-\frac{\sqrt3}{\sqrt{10}}\,e^5\wedge e^7\;.
\end{aligned}
\end{equation}
Note that the map $\alpha$ has no $\G$-symmetry.

With a similar trick as in Lemma~\ref{l2.1}, we can identify Clifford
multiplication on spinors with Cayley multiplication.  Recall that a
quotient~$M=G/H$ of compact Lie groups is equivariantly spin if and
only if the isotropy representation $\pi\colon H\to\SO(\frp)$ factors
over the spin group $\Spin(\frp)$.  In this case, the equivariant
spinor bundle $\Cal S\to M$ is constructed as the fibered product
$$
	\Cal S = G\times_\pix S\to M\;,
$$
where $\pix$ is the pull-back to $H$ by $\pi$ of the spin
representation of $\Spin(\frp)$ on the spinor module~$S$.  Since
Clifford multiplication $\frp\times S\to S$ is
$\Spin(\frp)$-equivariant, it is in particular $H$-equivariant, so
there is a fiber-wise Clifford multiplication $TM\times\Cal S\to\Cal
S$.

Note that if $\frp$ is odd-dimensional, Clifford multiplication with
vectors is uniquely defined only up to sign. To remove this ambiguity,
let
$$
	\omega =i^{\left[\frac{n+1}2\right]}e_1\cdots e_n
	\in\Cl(\frp)\otimes\C
$$
be the {\em complex Clifford volume element,\/} which
satisfies $\omega^2=1\in\Cl(\frp)$.  If $\frp$ is odd-dimensional,
then $\omega$ commutes with Clifford multiplication, and we require
that $\omega$ acts on $S$ as~$+1$.

The spinor module~$S$ of~$\Spin(7)$ is of dimension
$2^{\left[\frac72\right]}=8$. Because the smallest representations
of $\G$ are the trivial representation and the seven-dimensional
representation on the imaginary octonions $\I$, and because $\pix$ is
non-trivial, it is clear that there is a $\G$-equivariant
isomorphism $\Ca\cong\R\oplus\I\cong S$.

\begin{lem}\label{l2.2}
We identify $\frp\cong\I$ as in Lemma~\ref{l2.1}. With respect to a suitable
isometric, $\G$-equivariant identification $S\cong\Ca$ and a suitable
orientation of $\frp$, Clifford multiplication $\frp\times S\to S$
equals Cayley multiplication $\I\times\Ca\to\Ca$ from the right.
\end{lem}

In other words, the Clifford algebra $\Cl(\frp)\subset \End(S)
\cong\End_\R(\Ca)$ is generated by the endomorphisms $c_v$ given by
right multiplication with some element of $\frp\cong\I$.  Since $\Ca$
is not associative, we do not have the identity $c_v\cdot
c_w=c_{-v*w}$ in general.  For the same reason, right multiplication
on $\Ca$ does not commute with left multiplication, which agrees with
the fact that $S$ is an irreducible $\Cl(\frp)$-module.

\begin{proof} We fix a $\G$-equivariant orthogonal identification
$S\cong\Ca=\R\oplus\I$.  Clifford multiplication
$\frp\times(\R\oplus\I) \to(\R\oplus\I)$ splits into four components.
By Schur's lemma, the component $\frp\times\R\to\R$ vanishes.  Because
multiplication with a unit vector is an isometry on $S$, we have
$v\cdot 1=\pm v\in\I$ for $1\in\R\subset S$, and we choose the
identification $S\cong\R\oplus\I$ such that $v\cdot 1=v$.

Again by Schur's lemma, the component of $v\cdot s$ in $\R$
is $c\<v,s\>$ for some constant $c$. Since $v\cdot(v\cdot c)=-\norm
v^2\,c$, it is easy to see that $v\cdot v=-\norm v^2=v*v\in\R$
for $v\in\I$.

Finally, for orthogonal imaginary elements $v$, $w\in\I$ we must
have $\norm{v\cdot w}=\norm v\norm w=\norm{w*v}$, so $v\cdot w=\pm
w*v$.  To check that the correct sign is $+$, we calculate using
\eqref{2.3};
\begin{equation}\label{2.6}
	\omega\cdot s =\left(\cdots(s*e_7)*\cdots\right)*e_1
	=s\;.
\end{equation}
\end{proof}

\subsection{The spectra of some deformed Dirac operators}\label{s2.2}\hfill

We use the explicit formulas for Clifford multiplication and the
tangential part of the Lie bracket obtained in the previous section to
estimate the spectrum of the family of deformed odd signature
operators $B^{\lambda,3\lambda-1}$.
We take the orthonormal base~$e_1$, \dots, $e_7$ of~$\frp\cong\I$ as
in \eqref{2.2}. Let $c_i$ and $\hat c_i$ denote Clifford multiplication with
$e_i$ on the first and second factor of $\Lambda^\even\frp\cong
S\otimes S \cong\Ca\otimes\Ca$. Then the Clifford volume elements
$$
\omega=c_1\cdots c_7 \quad\text{ and }\quad \widehat
\omega=\hat c_1\cdots\hat c_7
$$
act as $1$ by \eqref{2.6}.

We extend $e_1$, \dots, $e_7$ to an orthonormal base $e_1$, \dots,
$e_{10}$ of $\frg$, and let~$c_{ijk}=\<[e_i,e_j]_\frp,e_k\>$, so for
example $c_{124}={\frac1{\sqrt5}}$ by Lemma~\ref{l2.1} and \eqref{2.3}.  We define
two symbols $\adX$ and $\adY\colon \frg\otimes \Lambda^\even \frp\to
\Lambda^\even \frp$ by
$$
	\aD_{\frp,i} = \aD_{\frp,e_i}
	=\frac1{ 4}\sum_{j,k=1}^mc_{ijk}\,c_jc_k
	\quad \text{ and }\quad
	\aDsl_{\frp,i} = \aDsl_{\frp,e_i}
	=\frac1{ 4}\sum_{j,k=1}^mc_{ijk}\,\hat c_j\hat c_k\;.
$$
Then $\pix=\adX|_{\frh}$ and $\piy=\adY|_{\frh}$ are the differentials
of the representations of $H$ on the two factors of $S\otimes S$ that
induce the bundle $\Lambda^\even TM\to M$.

We consider a family $D^\lambda$ of $G$-equivariant deformed Dirac
operators on $\Gamma(\Cal S)$ and a family $B^{\lambda,\mu}$ of
$G$-equivariant deformed odd signature operators on
$\Gamma(\Lambda^\even TM)$ as in \cite{g1}. Using Frobenius
reciprocity and the Peter-Weyl theorem, we will write
$$
\begin{aligned}
	\Gamma(\Cal S)
	&=\overline{\bigoplus_{\gamma\in\widehat G}
		V^\gamma\otimes\Hom_H(V^\gamma,S)}\\
	\quad \text{ and }\quad
	\Omega^\even(M)
	&=\overline{\bigoplus_{\gamma\in\widehat G}
		V^\gamma\otimes\Hom_H(V^\gamma,S\otimes S)}\;.
\end{aligned}
$$
Since $D^\lambda$ and $B^{\lambda,\mu}$ are $G$-equivariant, they
preserve these decompositions. Moreover, for each summand above, we
may write
$$
\begin{aligned}
	D^\lambda|_{V^\gamma\otimes\Hom_H(V^\gamma,S)}
	&=\id_{V^\gamma}\otimes\, ^\gamma\!D^\lambda\\
	\quad \text{ and }\quad
	B^{\lambda,\mu}|_{V^\gamma\otimes\Hom_H(V^\gamma,S\otimes S)}
	&=\id_{V^\gamma}\otimes\, ^\gamma\!B^{\lambda,\mu}\;.
\end{aligned}
$$
Let $\gamma_i$ denote the action of $\gamma^*_{e_i}$ on the dual of
the representation space $V^\gamma$.  With this notation, the
operators above take the form
\begin{equation}\label{2.7}
\begin{aligned}
	^\gamma\!D^\lambda
	&=\sum_{i=1}^7c_i\,\left(\gamma_i+\lambda\,\aD_{\frp,i}\right)\\
	\quad \text{ and }\quad
	 ^\gamma\!B^{\lambda,\mu}
	&=\sum_{i=1}^7c_i\,\left(\gamma_i+\lambda\,\aD_{\frp,i}+\mu\,
	\aDsl_{\frp,i}\right)\;.
\end{aligned}
\end{equation}
Note that $D=D^{\frac12}$ and $B=B^{{\frac12},{\frac12}}$ are
respectively the Dirac operator and the odd signature operator
associated to the Levi-Civita connection on $M$.  On the other hand,
$\Dsl=D^{\frac13}$ and $\Bsl=B^{{\frac13},0}$ are reductive operators
in the terminology of \cite{g1}, \cite{g2} and \cite{g3}.

We now consider the one-parameter family $B^{\lambda,3\lambda-1}$
for $\lambda\in\bigl[{\frac13},{\frac12}\bigr]$.  We write
\begin{equation}\label{2.8}
	^\gamma\!B^{\lambda,3\lambda-1}
	=\, ^\gamma\!\widetilde B
	+\mu\,\sum_{i=1}^7c_i\,\left({\frac13}\,\aD_{\frp,i} + \aDsl_{\frp,i}
	\right)
\end{equation}
for $\mu=3\lambda-1\in\left[0,{\frac12}\right]$.
Let us define
\begin{equation}\label{2.9}
B_0=\sum_{i=1}^7c_i\,\left({\frac13}\,\aD_{\frp,i} + \aDsl_{\frp,i} \right)\;.
\end{equation}

The square of $^\gamma\widetilde B$ has been computed in \cite{g1},
\cite{g2} as
\begin{equation}\label{2.10}
	^\gamma\!\widetilde B^2
	=\norm{\gamma+\rho_G}^2-c_H^{\piy}-\norm{\rho_H}^2\;.
\end{equation}
Here $\rho_H$ and $\rho_G$ are half sums of positive roots, and
$c_H^{\piy}$ is the Casimir operator of $H$ associated to the
representation $\piy$, taken with respect to the norm on $\frh$ that
is induced by the embedding~$\iota$ of \eqref{2.2} and a fixed Ad-invariant
scalar product on $\frg$.


Let $\frs\subset\frt$ be the Cartan subalgebras of $\frh\subset\frg$
spanned by $\iota_{12}$ and by $e_{12}$ and~$e_{34}$, respectively.
The weights of $\frh$ and $\frg$ are of the form
\begin{equation}\label{2.11}
ik\,\iota^*_{12}=\frac{ik}5\,(2e_{12}^*+e_{34}^*)\in i\frs^*
	\quad \text{ and }\quad
	ip\,e_{12}^*+iq\,e_{34}^*\in i\frt^*
\end{equation}
with~$k$, $p$, $q\in\Z$.
We will pick the Weyl chambers
\begin{equation}\label{2.12}
	P_H=\{\,it\,\iota^*_{12}\mid t\ge0\,\}\subset i\frs^*
	\quad \text{ and }\quad
	P_G=\{\,ix\,e_{12}^*+iy\,e_{34}^*\mid x\ge y\ge 0\,\}\;.
\end{equation}
With respect to these Weyl chambers, the dominant weights of $G$ and
$H$ are the weights in \eqref{2.11} with~$k\ge 0$ and~$p\ge q\ge 0$,
respectively.  Then we find
\begin{equation}\label{2.13}
	\rho_H=\frac i2\,\iota^*_{12}
	=\frac i{10}\,(2e_{12}^*+e_{34}^*)
	\quad \text{ and }\quad
	\rho_G=\frac i2\,(3e^*_{12}+e^*_{34})\;.
\end{equation}

Let $\gamma_{(p,q)}$ denote the irreducible $G$-representation with
highest weight $ip\,e_{12}^*+iq\,e_{34}$, where~$p\ge q\ge 0$ are
integers.  Let $\kappa_k$ denote the irreducible $H$-representation
with highest weight $ik\,\iota^*_{12}$, then the dimension of
$\kappa_k$ is $2k+1$.  We have seen above that the isotropy
representation $\pi$ on $\frp$ is isomorphic to $\kappa_3$, while
$\pix$ on $S\cong\R\oplus\I$ is isomorphic to $\kappa_0 \oplus
\kappa_3$.  We conclude that for $\gamma = \gamma_{(p,q)}$, we have
\begin{equation}\label{2.14}
	^\gamma\!\widetilde B^2
	=\begin{cases}
		\norm{\gamma_{(p,q)}+\rho_G}^2
			-\norm{\rho_H}^2
		=p^2+3p+q^2+q+\frac{49}{20}
			&\text{on~$\Hom_H(V^\gamma,S\otimes\R)$, and}\\
		\norm{\gamma_{(p,q)}+\rho_G}^2
			-\norm{\kappa_3+\rho_H}^2
		=p^2+3p+q^2+q+{\frac1{20}}
			&\text{on~$\Hom_H(V^\gamma,S\otimes\I)$.}
	\end{cases}
\end{equation}

We now calculate the spectral radii of the various components of the
operator $B_0$.  Since the operator $B_0$ evidently commutes with the
action of $\G$ on $S\otimes S$ by its definition in \eqref{2.9}, we
can restrict our attention to the $\G$-isotypic components of~$B_0$.
Let $\fru \subset \g$ be a Cartan subalgebra containing $\frs$.  We
introduce a basis of $i\fru^*\subset\fru^*\otimes_\R\C$ such that
$(1,0)$ and $(0,1)$ describe a long and a short root of $\g$
respectively, which belong to the closure of a fixed Weyl chamber in
$i\fru^*$.  In this basis, the dominant weights of $\g$ are given
precisely by pairs of non-negative integer coordinates.  Let
$\phy_{(a,b)}$ denote the irreducible $\G$-representation with highest
weight $(a,b)$ for $a$, $b\in\Z$ with $a$, $b\ge 0$. It is easy to
check that $\phy_{(0,1)}$ denotes the standard representation of $\g$
on $\I$, that $\phy_{(1,0)}$ is the adjoint representation, and that
$\phy_{(0,2)}$ is the 27-dimensional non-trivial part of the symmetric
product $S^2\I$.

Using the computer program LiE, we see that $S\otimes S$ splits
into $\G$- and $H$-isotypical components as:
\begin{equation}\label{2.15}
\begin{matrix}
	&\R\otimes\R
	&\cong_\G&\phy_{(0,0)}
	&\cong_H&\kappa_0\;,\\
\noalign{\vspace{2mm}}
	&\I\otimes\R
	&\cong_\G&\phy_{(0,1)}
	&\cong_H&\kappa_3\;,\\
\noalign{\vspace{2mm}}
	&\R\otimes\I
	&\cong_\G&\phy_{(0,1)}
	&\cong_H&\kappa_3\;,\\
\noalign{\vspace{3mm}}
	\text{and }\quad
	&\I\oplus\I
	&\cong_\G&\begin{matrix}
		\phy_{(0,0)}\oplus\phy_{(0,1)}\oplus\phy_{(1,0)}\\
		\oplus\phy_{(0,2)}\end{matrix}
	&\cong_H&\begin{matrix}
		\kappa_0\oplus\kappa_3\oplus(\kappa_1\oplus\kappa_5)\\
		\oplus(\kappa_2\oplus\kappa_4\oplus\kappa_6)\;.\end{matrix}\\
\end{matrix}
\end{equation}
Note that no two $\G$-representations involved have a common
isomorphic $H$-subrepresentation.

Lemmas~\ref{l2.1} and~\ref{l2.2} give us an explicit formula for $B_0$. Using a
computer program, we can calculate the eigenvalues of $B_0$ on each
$\G$-isotypical component $B_0^{(p,q)}$.  A basis for the trivial
component is given by
$$
1\otimes 1\quad\text{ and }\quad \frac1{\sqrt7}\,
\sum_{i=1}^7e_i\otimes e_i\;.
$$
With respect to this basis, one has
\begin{equation}\label{2.16}
B_0^{(0,0)} =\frac1{2\sqrt 5}\begin{pmatrix} 7&-3\,\sqrt 7\\-3\,\sqrt
	7&5\end{pmatrix}\;.
\end{equation}
In particular, the eigenvalues of $B_0^{(0,0)}$ are $\frac7{\sqrt5}$
and $-{\frac1{\sqrt5}}$.

The representation $\phy_{(0,1)}$ has multiplicity $3$ in
$\Lambda^\even\frp$.  We pick three vectors that equivariantly span
the isotypical component, and that correspond to~$e_1\in\I$, namely
$$
	e_1\otimes 1\;,\qquad
	1\otimes e_1\;,\quad\text{ and }\quad
	\frac1{\sqrt6}\,\sum_{i=2}^7e_i\otimes(e_1*e_i)\;.
$$
By equivariance, $B_0^{(0,1)}$ preserves the $3$-dimensional
subspace~$V\subset\Lambda^\even\frp$
spanned by these vectors.
We find
\begin{equation}\label{2.17}
	B_0^{(0,1)}|_V
	=\frac1{2\sqrt 5}
		\begin{pmatrix} -1&3&3\,\sqrt6\\
			3&7&\sqrt6\\
			3\sqrt6&\sqrt6&-4\end{pmatrix}\;.
\end{equation}
The eigenvalues of $B_0^{(0,1)}$ are readily computed to be
${\frac1{\sqrt5}}$ and $\pm\sqrt5$.

Finally, the $\G$-isotypical components isomorphic to $\phy_{(1,0)}$
and $\phy_{(0,2)}$ both have multiplicity $1$, and we have
\begin{equation}\label{2.18}
	B_0^{(1,0)}={\frac1{\sqrt5}}\quad\text{ and }\quad
	B_0^{(0,2)}=-{\frac1{\sqrt5}}\;.
\end{equation}
The calculations above lead to the following proposition.

\begin{prop}\label{p2.3}
The operator~$D^\lambda$ has no kernel for $\lambda \in
\bigl[{\frac13}, {\frac12} \bigr]$.  For $\lambda\in \bigl[{\frac13},
{\frac12}\bigr]$ and $\gamma\in \widehat G$, the operator $^\gamma\!
B^{\lambda, 3\lambda-1}$ has a non-zero kernel only if $\lambda =
{\frac12}$ and $\gamma = \gamma_{(0,0)}$ is the trivial
representation.  For $\gamma = \gamma_{(0,0)}$, the operator $^\gamma
\!B^{\lambda,3\lambda-1}$ has a positive and a negative eigenvalue if
$\lambda \in \bigl[{\frac13}, {\frac12}\bigr)$, and only the negative
eigenvalue vanishes at $\lambda={\frac12}$. \end{prop}

\begin{proof} 
The claim about $D^\lambda$ follows from the proof of Lemma 4.6 in
\cite{g2} (see also Bemerkung 1.20 in \cite{g1}).

Let us now check that for no non-trivial representation $\gamma$ of
$G$ and no $\lambda \in \bigl[{\frac13}, {\frac12}\bigr]$, the
operator $^\gamma\!B^{\lambda, 3\lambda-1}$ can have a kernel.  This
is because by \eqref{2.14}, all eigenvalues of $^\gamma\!\widetilde B$
belong to $\R\setminus \bigl(-\frac9{2\sqrt 5}, \frac9{2\sqrt
5}\bigr)$, where $\pm\frac9{2\sqrt 5}$ is attained on
$\Hom_H(V^\gamma, S\otimes\frp)$ for $\gamma = \gamma_{(1,0)}$.  On
the other hand, the spectral radius of $\mu\,B_0$ is
$\frac{7\mu}{\sqrt5}$, which is smaller than $\frac9{2\sqrt5}$ for
$\mu \in\bigl[0,{\frac12} \bigr]$.

Now, consider the operator $^\gamma\!\widetilde B$ for the trivial
representation $\gamma = \gamma_{(0,0)}$.  Clearly,
$\Hom_H(V^\gamma,S\otimes S)$ is isomorphic to the trivial
$\G$-isotypical component of $\Lambda^\even\frp$.  Another machine
computation shows that in the basis of \eqref{2.16}, the operator
$^\gamma\!\widetilde B$ takes the form
$$
	^\gamma\!\widetilde B
	=\frac1{2\sqrt5}\,\begin{pmatrix}7&\\&-1\end{pmatrix}\;.
$$
The eigenvalues of the operator
$$
	^\gamma\!B^{\frac{\mu+1}3,\mu}
	=^\gamma\!\widetilde B+\mu\,B_0
	=\frac1{2\sqrt5}
		\,\begin{pmatrix}7+7\mu&-3\mu\,\sqrt7\\
			-3\mu\,\sqrt7&5\mu-1\end{pmatrix}
$$
are precisely $\frac1{2\sqrt5}\,\bigl(6\mu + 3\pm\sqrt{64 \mu^2 + 8\mu
+ 16}\bigr)$.  Since $6\mu+3 < \sqrt{64\mu^2 + 8\mu + 16}$ except at
$\lambda={\frac12}$ where one gets equality, the claims in the
proposition follow.  \end{proof}

\subsection{Computing the $\eta$-invariants.}\label{s2.3}\hfill

Next we compute the $\eta$-invariants $\eta(B)$ and $\eta(D)$ for the
Dirac operators considered in the previous subsection, up to a local
correction term.  We will use the formula of \cite{g3}.

We fix Weyl chambers $P_G$ and $P_H$ as in \eqref{2.12}. Then $\rho_G$ and
$\rho_H$ are given by \eqref{2.13}.  The choices of $P_G$ and $P_H$ also
determine orientations on~$\frg/\frt$ and $\frh/\frs$.  If $\alpha_1$,
\dots, $\alpha_l\in i\frt^*$ are the positive roots of $\frg$ with
respect to $P_G$, then we can choose a complex structure on
$\frg/\frt$ and a complex basis $z_1$, \dots, $z_n$ such that
$\ad|_{\frt\times(\frg/\frt)}$ takes the form
$$
	\ad_X=\begin{pmatrix} \alpha_1(X)\\&\ddots\\&&\alpha_l(X)
	\end{pmatrix} \quad \text{ for all~$X\in\frt$.}
$$
Then we declare the real basis $z_1$, $i\,z_1$, $z_2$, \dots, $i\,z_l$
to be positively oriented.

Having fixed orientations on $\frp=\frg/\frh$ by a choice of
an orthonormal base in \eqref{2.4}
and orientations on $\frg/\frt$ and $\frh/\frs$ as above,
there is a unique orientation on $\frt/\frs$
such that the orientations on
$$
	\frg/\frs
	\cong\frp\oplus(\frh/\frs)
	\cong(\frg/\frt)\oplus(\frt/\frs)
$$
agree. Let $E\in \frt/\frs \cong \frs^\perp \subset \frt$ be the
positive unit vector, and let $\delta\in i\frt^*$ be the unique weight
such that
	$$-i\delta(E)>0 \quad\text{ and }\quad \delta(X)\in2\pi i\Z
	\iff e^X\in S$$
for all $X\in\frt$. Then one can check that
\begin{equation}\label{2.19}
	E={\frac1{\sqrt5}}\,\left(e_{12}-2e_{34}\right)
		\quad\text{ and }\quad
	\delta=i\,\left(e_{12}^*-2e_{34}^*\right)
\end{equation}
are compatible with the orientations fixed above.

Let $\Dsl^{(k)}$ be the reductive Dirac operator acting on
$\Gamma(\Cal S\otimes V^\kappa M)$, where $\kappa$ is the
$H$-representation with highest weight $\kappa_k=ik\,\iota_{12}^*$.
Then we note that
$$
\Dsl=\Dsl^{(0)} \quad\text{ and }\quad \Bsl=\Dsl^{(0)}\oplus\Dsl^{(3)}\;.
$$
We have to find the unique weights $\alpha_k\in i\frt^*$ of $\frg$ such that
$$
	\alpha_k|_\frs=ik\,\iota_{12}^*+\rho_H
	\quad\text{ and }\quad
	-i(\alpha_k-\delta)(E)<0\le-i\alpha_k(E)\;.
$$
By \eqref{2.19} we have
$$
	\alpha_0=i\left({\frac12}\,e_{12}^*-{\frac12}\,e_{34}^*\right)
	\quad\text{ and }\quad
	\alpha_3=i\left({\frac32}\,e_{12}^*+{\frac12}\,e_{34}^*\right)\;.
$$
Note that $\alpha_0(E)$, $\alpha_3(E)\ne 0$.

Let $\Delta_+=\{ie_{12}^*+ie_{34}^*,ie_{12}^*-ie_{34}^*,ie_{12}^*,ie_{34}^*\}$
denote the set of positive roots with respect to $P_G$.
Let $\adach$ denote the map, $z\mapsto\frac z{2\sinh(z/2)}$.

We also need some equivariant characteristic differential forms.  Note
that we will eventually evaluate these forms only at~$X=0$, so that we
may actually forget the equivariant formalism in a moment.  Let
$\adach_X(M,\nabla)$ be the total equivariant $\adach$-form, and
$\ldach_X(M,\nabla) = 2\, \adach_X(M,\nabla) \wedge \ch_X(\Cal
S,\nabla)$ be a rescaled equivariant $L$-form, both taken with respect
to a connection $\nabla$ on $TM$ and the induced connection on $\Cal
S$.
If $p_k=p_k(M,\nabla)\in \Omega^{4k}_\frg(M)$ denotes the $k$-th
equivariant Pontrjagin form of $M$, then
$$
	\Adach_X(M,\nabla)
	=1-\frac{p_1}{24}+\frac{7\,p_1^2-4\,p_2}{2^7\,3^2\,5}+\dots
	\quad\text{ and }\quad
	\Ldach_X(M,\nabla)
	=16+\frac{4\,p_1}{3}+\frac{7\,p_2-p_1^2}{45}+\dots
$$
So in particular,
\begin{equation}\label{2.20}
	\Adach_X(M,\nabla)+\frac{\Ldach_X(M,\nabla)}{2^5\,7}
	=\frac{15}{14}-\frac{p_1(M,\nabla)}{28}
	+\frac{p_1(M,\nabla)^2}{2^7\,7}+\dots
\end{equation}
For different connections $\nabla$, $\nabla'$ let
$\adachsl(M,\nabla,\nabla')$ and $\ldachsl(M, \nabla,\nabla')
\in\Omega^*_\frg(M)/d_\frg \Omega^*_\frg(M)$ denote the corresponding
equivariant Chern-Simons classes with
$$
	d\Adachsl(M,\nabla,\nabla')=\Adach_X(M,\nabla')-\Adach_X(M,\nabla)
	\quad\text{ and }\quad
	d\Ldachsl(M,\nabla,\nabla')=\Ldach_X(M,\nabla')-\Ldach_X(M,\nabla)\;.
$$
We will work with the reductive connection $\nabla^0$ and the Levi-Civita
connection $\nabla\LC$ on $TM$.

We can now compute the eta-invariants of $D$ and $B$ using the formula
for infinitesimal equivariant $\eta$-invariants computed in \cite{g3}.

\begin{thm}\label{t2.4} The $\eta$-invariants of $D$ and $B$ are the values at
$X=0 \in\frt$ of the following:
\begin{align*}
\begin{split}
	\eta_X(D)
	&=2\,\sum_{w\in W_G}\frac{{\rm sign}(w)}{\delta(wX)}\,\Biggl(
		\prod_{\beta\in\Delta_+}\Adach\left(\beta(wX)\right)
		\cdot\Adach\left(\delta(wX)\right)
			\,e^{\left(\alpha_0-\frac\delta2\right)(wX)}\\
	&\hspace{3.5cm}
		-\prod_{\beta\in\Delta_+}\Adach\left(\beta(wX|_\frs) \right)
		\cdot e^{\rho_H(wX|_\frs)}\Biggr)
		\cdot\prod_{\beta\in\Delta_+}\frac1{\beta(X)}\\
	&\qquad
		+\int_M\Adachsl_X\left(TM,\nabla^0,\nabla\LC \right)\;,
\end{split}\tag1\\
\noalign{\medskip}
\begin{split}
	\eta_X(B)
	&=2\,\sum_{w\in W_G}\frac{{\rm sign}(w)}{\delta(wX)}\,\Biggl(
		\prod_{\beta\in\Delta_+}\Adach\left(\beta(wX)\right)
		\cdot\Adach\left(\delta(wX)\right)
			\,\left(e^{\left(\alpha_0-\frac\delta2\right)(wX)}
				+e^{\left(\alpha_3-\frac\delta2\right)(wX)}
			\right)\\
	&\hspace{3.5cm}
		-\prod_{\beta\in\Delta_+}\Adach\left(\beta(wX|_\frs)\right)
		\cdot\left(e^{\rho_H(wX|_\frs)}
			+e^{(\kappa_3+\rho_H)(wX|_\frs)}\right)
		\Biggr) \cdot\prod_{\beta\in\Delta_+}\frac1{\beta(X)}\\
	&\qquad
		+1+\int_M\Ldachsl_X\left(TM,\nabla^0,\nabla\LC \right)\;.
\end{split}\tag2
\end{align*}
\end{thm}

\begin{proof} This follows immediately from \cite{g3},
Theorem~2.33 and Corollary~2.34,
and from Proposition~\ref{p2.3}. \end{proof}

A machine calculation now gives numerical values up to the local
correction term.

\begin{cor}\label{c2.5} We have the formulae
$$
\begin{aligned}
	\eta(D)
	&=-\frac{12923}{2\;3^2\;5^6}
		+\int_M\Adachsl_X\left(TM,\nabla^0,\nabla\LC \right)\;,\\
		\text{and }\quad
	\eta(B)
	&=1-\frac{12923}{2\;3^2\;5^6}-\frac{277961}{2\;3^2\;5^6}
		+\int_M\Ldachsl_X\left(TM,\nabla^0,\nabla\LC \right)\\
	&=-\frac{4817}{3^2\;5^6}
		+\int_M\Ldachsl_X\left(TM,\nabla^0,\nabla\LC \right)\;.
\end{aligned}
$$
\end{cor}

\begin{rem}\label{r2.6}
One might be tempted to conjecture that $\eta(\Dsl)$ has the
value~$-\frac{12923}{2\;3^2\;5^6}$ (from above) because $\Dsl^2$
involves the Laplacian on $\Cal S$ with respect to the reductive
connection $\nabla^0$.  However, the equivariant $\eta$-invariant
$\eta_G(\Dsl)$ has been calculated in \cite{g1} and in particular
$\eta(\Dsl) = \frac{207479}{2^5\;3^2\;5^6} \ne -
\frac{12923}{2\;3^2\;5^6} = -\frac{206768}{2^5\;3^2\;5^6}$.
\end{rem}

\subsection{Equivariant $\eta$-invariants and the Eells-Kuiper
invariant.}\label{s2.4}\hfill

We compute the Eells-Kuiper invariant of $M=\SO(5)/\SO(3)$ using
Donnelly's formula \cite{donnelly}; see also \cite{ks} which involves
the non-equivariant $\eta$-invariants of the Dirac operator $D$ and
the signature operator $B$ on $M$.  Using the methods of \cite{g3}, we
determine $\eta(B)$ and $\eta(D)$ from their equivariant counterparts
computed in Theorem~\ref{t2.4} above for those group elements that act
freely.

Recall that the Eells-Kuiper invariant is defined as (see Section \ref{s2.1})
$$
	\ek(M)
	=\frac{\eta(B)}{2^5\,7}+\frac{\eta(D)+h(D)}2
		-\frac1{2^7\,7}\,\int_Mp_1(M,\nabla\LC)\wedge h(M,\nabla\LC)
	\qquad\in\Q/\Z\;,
$$
(see \eqref{EK}). Note that the form $h(M,\nabla\LC)$ exists because
$H^4(M,\R)=0$, and is unique up to exact forms because $H^3(M,\R)=0$.
Moreover, we may choose $h(M,\nabla\LC)$ to be $G$-invariant.

As above, let $\adachsl(M,\nabla,\nabla') \in\Omega^*(M)/d
\Omega^*(M)$ denote the Chern-Simons class that interpolates between
the $\adach$-forms constructed from two connections $\nabla$ and
$\nabla'$.  If $h=h(\nabla) \in\Omega^*(M)/d \Omega^*(M)$ is a class
such that $dh=p_1(M,\nabla)$ is the first Pontrjagin form of $TM$,
then the class $h(\nabla') = h(\nabla) + \tilde p_1(M,\nabla,\nabla')$
satisfies $dh(\nabla') = p_1(M,\nabla')$. In particular
\begin{equation}\label{2.21}
\begin{aligned}
	\frac1{ 2^7\,7}\int_M\biggl(p_1(M,\nabla')\,h(M, &\nabla')
		-p_1(M,\nabla)\,h(M,\nabla)\biggr)\\
	&=\frac1{ 2^7\,7}\int_M\left(\tilde p_1(M,\nabla,\nabla')\,
	p_1(M, \nabla)
		+p_1(M,\nabla')\,\tilde p_1(M,\nabla,\nabla')\right)\\
	&=\int_M\left(\Adachsl\left(M,\nabla,\nabla'\right)
		+\frac1{2^5\,7}\,\Ldachsl\left(M,\nabla,\nabla'\right)
	\right)
\end{aligned}
\end{equation}
by \eqref{2.20}.  As an immediate consequence of Corollary~\ref{c2.5},
\eqref{EK} and
\eqref{2.21}, we get
\begin{equation}\label{2.21a}
\begin{aligned}
	\ek(M)
	&=-\frac{12923}{2\;3^2\;5^6}-\frac{4817}{2^5\;7\cdot 3^2\;5^6}
		+\int_M\biggl(\Adachsl\left(M,\nabla^0,\nabla\LC\right)
		+\frac1{2^5\;7}\,\Ldachsl\left(M,\nabla^0,\nabla\LC
	\right) \biggr)\\
	&\hspace{2cm}
		-\frac1{ 2^7\,7}\int_Mp_1(M,\nabla\LC)\,h(M,\nabla\LC)\\
	&=-\frac{16189}{2^5\,5^5\,7}
		-\frac1{ 2^7\,7}\int_Mp_1(M,\nabla^0)\,h(M,\nabla^0)\;.\\
\end{aligned}
\end{equation}
Here we have used that~$D$ has no kernel by Proposition~\ref{p2.3}, so
$h(D)=0$.

\subsection{Computing the Eells-Kuiper invariant.}\label{s2.5}\hfill

It remains to evaluate the integral of the secondary class
$p_1(M,\nabla^0)\,h(M,\nabla^0)$ over $M$. This is again done with
the help of the results of Section \ref{s2.1}.

Let $V$, $W$ be vector fields on $M$.  Then there exist
$H$-equivariant functions $\widehat V$, $\widehat W\colon G\to\frp$ such that
$$
	V(gH)=\bigl[g,\widehat V(g)\bigr]
	\quad\text{ and }\quad
	W(gH)=\bigl[g,\widehat W(g)\bigr]\in TM=G\times_\pi\frp\;.
$$
The reductive connection $\nabla^0$ and its curvature $R^0$ satisfy
$$
	\widehat{\nabla^0_VW}=\widehat V\left(\widehat W\right)
	\quad\text{ and }\quad
	\widehat{R^0_{V,W}}=-\pi_{*\left[\widehat V,\widehat W\right]_\frh}\;.
$$
Because $p_1(M,\nabla^0)$ is $G$-invariant, it must be given by an
$H$-invariant $\hat p_1(M,\nabla^0)\in\Lambda^4\frp^*$.  Then $\hat
p_1(M,\nabla^0)$ is in fact $\G$-invariant by \eqref{2.15}, and hence,
it must be a multiple the Poincar\'e dual~$\lambda_4$ of the three
form~$\lambda_3$ where
\begin{equation}\label{2.22}
\begin{aligned}
	\lambda_3
	&=\<\punkt*_\I\punkt,\punkt\>
	=\sum_{i=1}^7e^i\wedge e^{i+1}\wedge e^{i+3}\;,\\
		\text{so }\quad
	\lambda_4
	&=\sum_{i=1}^7e^i\wedge e^{i+1}\wedge e^{i+2}\wedge e^{i+5}\;.
\end{aligned}
\end{equation}
It is thus sufficient to compute $\hat p_1(M, \nabla^0)
(e_2,e_4,e_5,e_6)$. Using \eqref{2.5}, one can check that $[e_2,e_5]$,
$[e_2,e_6]$, $[e_4,e_5]$, $[e_4,e_6]\in\frp$, and
$$
\begin{aligned}
	\hat p_1(M,\nabla^0)(e_2,e_4,e_5,e_6)
	&=-\frac1{8\pi^2}\,\tr\left(\widehat R^0)^2\right)(e_2,e_4,e_5,e_6)
	=-\frac1{4\pi^2}\,\tr\left(\pi_{*[e_2,e_4]_\frh}\,
	\pi_{*[e_5,e_6]_\frh}\right)\\
	&=\frac3{20\pi^2}\,\tr\left(\pi_{*f_1}^2\right)
	=-\frac{21}{25\pi^2}\;.
\end{aligned}
$$
This implies that
\begin{equation}\label{2.23}
	\hat p_1(M,\nabla^0)
	=-\hat p_1(M,\nabla^0)(e_2,e_4,e_5,e_6)\,\lambda_4
	=\frac{21}{25\pi^2}\,\sum_{i=1}^7e^i\wedge e^{i+1}
		\wedge e^{i+2}\wedge e^{i+5}\;.
\end{equation}
Clearly, \eqref{2.23} gives an $H$-invariant element of $\Lambda^4\frp^*$.
By equivariance, we can write $p_1(M, \nabla^0) = dh(M,\nabla^0)$ for
some $G$-invariant form $h(M,\nabla^0)$, which is again given by an
$H$-invariant $\widehat{h(M,\nabla^0)} \in \Lambda^3\frp^*$.  Thus,
$h(M,\nabla^0)$ must be a multiple of the form $\lambda_3$ of
\eqref{2.22}. By Cartan's formula for the exterior derivative and Lemma~\ref{l2.1},
we have
$$
\begin{aligned}
	\widehat{d\lambda_3}\left(\widehat V_0,\dots,\widehat V_3\right)
	&=-\widehat\lambda_3\left(\bigl[\widehat V_0,\widehat V_1\bigr],\widehat V_2,
		\widehat V_3\right)
	+\widehat\lambda_3\left(\bigl[\widehat V_0,\widehat V_2\bigr],\widehat V_1,
		\widehat V_3\right)\mp\dots
	-\widehat\lambda_3\left(\bigl[\widehat V_2,\widehat V_3\bigr],\widehat V_0,
		\widehat V_1\right)\\
	&=\frac6{\sqrt 5}\,\sum_{i=1}^7e^i\wedge e^{i+1}
		\wedge e^{i+2}\wedge e^{i+5}\;,
\end{aligned}
$$
so $\hat h(M,\nabla^0) = -\frac7{10 \sqrt 5\, \pi^2}\, \lambda_3$.
In particular, $p_1(M,\nabla^0)\,h(M,\nabla^0)$ is given by
$$
	\frac{7}{10\sqrt 5}\,\lambda_3\,\frac{21}{25\pi^2}
		\,\lambda_4
	=\frac{3\;7^3}{2\;5^{\frac72}\,\pi^4}\,e^1\wedge\dots\wedge
		e^7\;.
$$
Because $\Vol \left(\SO(3) \right) = 8\,\pi^2$,
$\Vol\left(\SO(5) \right) = \frac{2^7\,\pi^6}3$,
$\Vol(H) = 5^{\frac32}\, \Vol\left( \SO(3) \right)$, the volume of $M$
is~$\frac{16\pi^4}{3\;5^{\frac32}}$.  We can now calculate the last
contribution to $\ek(M)$ as
\begin{equation}\label{2.24}
	-\frac1{2^7\,7}\,\int_Mp_1(M,\nabla^0)\,h(M,\nabla^0)
	=-\frac1{2^7\;7\vphantom{5^{\frac72}}}
		\cdot\frac{3\;7^3}{2\;5^{\frac72}\,\pi^4}
		\cdot\frac{16\pi^4}{3\;5^{\frac32}}
	=-\frac{49}{50\,000}\;.
\end{equation}

Together with~\eqref{2.21a},
this completes the proof of Theorem~\ref{MainThm}. \hfill $\Box$

\section{Vector fields on 2-connected 7-manifolds}\label{s3}

In this section we prove some general results about smooth, oriented,
$2$-connected $7$-manifolds. In particular, we determine the maximal
number of independent vector fields on such a manifold in terms of the
first spin characteristic class. We will describe these results after
some general remarks about the homotopy type of $2$-connected
$7$-manifolds.

Let $M$ be a closed, oriented, $2$-connected $7$-manifold. From the structure
of $H^*(M)$, we see that $M$ is homotopy equivalent to a $CW$-complex
with cells in dimension $0,3,4$ and $7$. Furthermore, it follows from
Poincar\'{e} duality that the number of cells in dimension $3$ equals
the number of cells in dimension $4$, and that there is a unique cell
in dimensions $0$ and $7$. Let $M_k$ denote the $k$-skeleton of
$M$. We will denote by $M_k/M_{k-1}$ the space obtained by pinching
off the $(k-1)$-skeleton from $M_k$. Hence, $M_k/M_{k-1}$ is
equivalent to a one-point union of $k$-spheres.

\begin{prop}\label{p3.1} The following composite map is null homotopic
$$\begin{CD}
	\sph^6 @>>> M_4 @>>> M_4/M_3,
\end{CD}
$$
where the first map is the attaching map for the $7$-cell and the next
map is the pinch map.
\end{prop}
\begin{proof}
Since $M$ is $2$-connected, it admits a spin structure. The argument
given on page 32 of \cite{mm} shows that the above map is trivial for
any seven dimensional spin manifold.
\end{proof}

Any oriented, $2$-connected manifold admits a unique compatible spin
structure. Let $\beta \in H^4(M)$ be the first spin characteristic
class. The class $\beta$ is related to the first Pontrjagin class by
the relation $2\beta = p_1$. The relation to the Stiefel Whitney
classes is given by $\beta \equiv w_4\, (\!\!\!\mod 2)$. We now recall
the definition of the Wu classes.

\begin{defn}\label{d3.2}
For an $n$-manifold $M$, we define the Wu classes $\wu_i \in
H^i(M,\f_2)$ by the property $\wu_i\cup x = Sq^i(x)$ for all $x \in
H^{n-i}(M,\f_2)$, where $Sq^i$ denotes the $i$-th Steenrod operation.
\end{defn}
Let $\wu_t = 1 + \wu_1t + \wu_2t^2 + \ldots$, be the total Wu class. Then
the total Wu class is related to the total Stiefel Whitney class by
the relation $\wu_t \cup SW_t = 1$.

\begin{prop}\label{p3.3}
If $M$ is a $2$-connected $7$-manifold, then all of its Stiefel
Whitney classes are trivial.
\end{prop}
\begin{proof}
By the relation between the Stiefel Whitney classes and the Wu
classes, it is sufficient to show that the total Wu class is
trivial. Since the Steenrod operations are unstable cohomology
operations, it follows that $\wu_i = 0$ for $i>3$. By the sparseness of
$H^*(M,\f_2)$, it follows that $\wu_1 = \wu_2 = 0$. Finally, one uses
the Adem relation $Sq^3 = Sq^1Sq^2$ to see that $\wu_3 = 0$.
\end{proof}

\begin{cor}\label{c3.4}
The class $\beta \in H^4(M)$ is $2$-divisible. In other words, there
is some (not necessarily unique) class $\gamma$ such that $2\gamma =
\beta$.
\end{cor}
\begin{proof}
We know that $\beta \equiv w_4\, (\!\!\!\mod 2)$. Since $w_4 = 0$, $\beta
\equiv 0\,(\!\!\!\mod 2)$ which says that $\beta$ is $2$ divisible.
\end{proof}

We now proceed to use the above facts to study the vector fields on
$M$.
\begin{thm}\label{t3.5}
Any smooth, oriented, $2$-connected $7$-manifold is parallelizable if
and only if $\beta = 0$.
\end{thm}
\begin{proof}
Let $f : M \rightarrow B_{\Spin(7)}$ be the map that classifies the
tangent bundle of $M$. We would like to show that $f$ is homotopic to
the trivial map. The obstructions to constructing the null homotopy
lie in the groups $H^{i+1}(M,\pi_i(\Spin(7)))$. The primary obstruction
lies in $H^4(M,\pi_3(\Spin(7))) = H^4(M)$ and is none other than the
class $\beta$ which we assumed to be trivial. Since $\pi_6(\Spin(7)) =
0$, there are no further obstructions to constructing the null homotopy.
\end{proof}

\begin{thm}\label{t3.6}
Any non-parallelizable smooth, oriented, $2$-connected $7$-manifold
admits exactly $4$ independent vector fields. In other words, the
structure group of $M$ may be reduced to $\Spin(3)$ and no further.
\end{thm}
\begin{proof}
Let us first see that the structure group cannot be reduced further
than $\Spin(3)$. Suppose we could reduce the structure group to
$\Spin(2)$. Since $\Spin(2) = S^1$, any map from $M$ to $B_{\Spin(2)}$
is classified by $H^2(M)$. But $M$ is $2$-connected so any such map is
trivial. Hence, a reduction of the structure group to $\Spin(2)$ would
mean that the manifold is parallelizable.
 
It remains to show that we can always reduce the structure group to
$\Spin(3)$. This corresponds to lifting the map $f : M \rightarrow
B_{\Spin(7)}$ to the space $B_{\Spin(3)}$. The obstructions to
constructing this lift lie in the groups $H^{i+1}(M, \pi_i
(\Spin(7)/\Spin(3)))$. The primary obstruction lies in the group
$H^4(M,\f_2)$ and is, by naturality, the element $\beta\, (\!\!\!\mod
2)$. Since $\beta$ is $2$ divisible, this element is zero. Hence, we
can construct the lift on $M_4$, the $4$-skeleton of $M$. One now has
the following commutative diagram:
$$
\xymatrix{
 \sph^6 \ar@{->}[r]^{g\phantom{abcdefgh}} \ar[d] & {\Spin(7)/\Spin(3)} \ar[d] \\
 M_4 \ar[r] \ar[d] & B_{\Spin(3)} \ar[d] \\
 M \ar[r]_{f\phantom{abcd}}          & B_{\Spin(7)} 
}
$$
where the vertical maps on the left form a cofibration sequence and
those on the right form a fibration sequence. Let $g$ denote the map
$\sph^6 \rightarrow \Spin(7)/\Spin(3)$. In order to complete the lift to
all of $M$, we require that the map $g$ is null homotopic. Consider
the composite $\sph^6 \rightarrow \Spin(7)/\Spin(3) \rightarrow
B_{\Spin(3)}$. Since $B_{\Spin(3)}$ is $3$-connected, this composite
factors through the map $\sph^6 \rightarrow M_4/M_3$, which we know is
trivial. Hence, we know that $g$ lifts to $\Spin(7)$. But $\pi_6
(\Spin(7)) = 0$, and so $g$ is null homotopic.
\end{proof}

\begin{rem}\label{r3.7} From the above discussion, the Berger space admits
precisely four independent vector fields. \end{rem}

\bibliographystyle{alpha}

\begin{thebibliography}{Abc}

\bibitem[APS75]{aps} M. F. Atiyah, V. K. Patodi, I. M. Singer,
Spectral asymmetry and Riemannian geometry. I,
\textit{Math.\ Proc.\ Cambridge\ Philos.\ Soc.,} 77 (1975), 43--69

\bibitem[Ber61]{berger} M.\ Berger, 
Les vari\'et\'es riemanniennes homog\`enes normales simplement
connexes \`a courbure strictement positive,
\textit{Ann.\ Scuola Norm.\ Sup.\ Pisa (3),} 15 (1961), 179--246.

\bibitem[CLL99]{lie} A.\ M.\ Cohen, M.\ v.\ Leeuwen, B.\ Lisser,
\textit{LiE---a software package for Lie group
computations, version~2.2,} 1999.

\bibitem[CE00]{ce} D.\ Crowley, C.\ Escher,
Classification of $\sph^3$-bundles over $\sph^4$,
preprint (2000), {\tt arXiv:math.AT/0004147}.

\bibitem[Don75]{donnelly} H.\ Donnelly,
Spectral geometry and invariants from differential topology,
\textit{Bull.\ London Math.\ Soc.,} 7 (1975), 147--150.

\bibitem[DR81]{dr} A.\ Derdzinski, A.\ Rigas,
Unflat connections in 3-sphere bundles over $\sph^4$,
\textit{Trans.\ of the AMS,} 265 (1981), 485--493.

\bibitem[EK62]{ek} J.\ Eells, Jr., N.\ H.\ Kuiper,
An invariant for certain smooth manifolds,
\textit{Ann.\ Mat.\ Pura Appl., (4)} 60 (1962), 93--110.

\bibitem[Go97]{g1}S.\ Goette,
\textit{\"Aquivariante $\eta$-Invarianten homogener R\"aume,}
Shaker, Aachen (1997).

\bibitem[Go99]{g2} S.\ Goette,
Equivariant $\eta$-invariants on homogeneous spaces,
\textit{Math.\ Z.,} 232 (1999), 1--42.


\bibitem[Go02]{g3} S.\ Goette,
$\eta$-invariants of homogeneous spaces,
preprint (2002), {\tt arXiv:math.DG/0203269}.

\bibitem[GM74]{gm} D.\ Gromoll, W.\ Meyer,
An exotic sphere with nonnegative sectional curvature,
\textit{Ann.\ of Math.,} 100 (1974), 401--406.

\bibitem[GZ00]{gz} K.\ Grove, W.\ Ziller,
Curvature and symmetry of Milnor spheres,
\textit{Ann.\ of Math.,} 152 (2000), 331--367.

\bibitem[He60]{hermann} R. Hermann,
A sufficient condition that a mapping of Riemannian manifolds
be a fibre bundle,
\textit{Proc. Am. Math. Soc.,} 11 (1960), 236--242.

\bibitem[Kla88]{klaus} S.\ Klaus, 
\textit{Einfach-Zusammenh\"angende Kompakte Homogene R\"aume bis zur
Dimension Neun,} Diploma thesis, University of Mainz, 1988.

\bibitem[Kos99]{kostant} B.\ Kostant,
A cubic Dirac operator and the emergence of Euler number multiplets
of representations for equal rank subgroups,
\textit{Duke Math.\ J.,} 100 (1999), 447--501.

\bibitem[KS88]{ks} M.\ Kreck, S.\ Stolz,
A diffeomorphism classification of 7-dimensional homogeneous Einstein manifolds
with $\SU(3)\times\SU(2)\times\U(1)$-symmetry,
\textit{Ann.\ of Math.,} 127 (1988), 373--388.

\bibitem[KiSh01]{kish} N.\ Kitchloo, K.\ Shankar,
On complexes equivalent to $\sph^3$-bundles over $\sph^4$,
\textit{Int.\ Math.\ Research Notices,} 8 (2001), 381--394.

\bibitem [MM79]{mm} I.\ Madsen, R.\ J.\ Milgram, 
\textit{The classifying spaces for surgery and cobordism of
manifolds,} Princeton University Press, 1979.

\bibitem[On66]{onishchik} A.\ L.\ Onishchik, Transitive compact
transformation groups, \textit{Amer.\ Math.\ Soc.\ Transl.}
\textbf{55} (1966), 153--194.


\bibitem[Sle87]{sleb} S.\ Slebarski,
The Dirac operator on homogeneous spaces
and representations of reductive Lie groups I,
\textit{Am.\ J.\ Math.} 109 (1987), 283--302.

\bibitem[Str94]{straume} E.\ Straume, Compact differentiable
transformation groups on exotic spheres, \textit{Math.\ Ann.}
\textbf{299} (1994), 355--389.

\bibitem[Wo68]{wolf} J.\ Wolf,
The geometry and structure of isotropy irreducible homogeneous spaces,
\textit{Acta Math.,} 120 (1968), 59--148; correction in \textit{Acta
Math.,} 152 (1984), 141--142.


\bibitem[Zil01]{ziller2} W.\ Ziller,
\textit{Homogeneous spaces and positive curvature,} Lecture notes, 2001.


\end{thebibliography}

\end{document}